\documentclass[12pt]{amsart}

\usepackage{longtable}
\usepackage{url}
\usepackage{enumerate}
\usepackage{epsfig}
 \usepackage{ifthen}
 \newcommand{\mymarginpar}[1]{%
    \marginpar{\ifthenelse{\isodd{\arabic{page}}}{\flushleft #1}{\flushright #1}}}

 \renewcommand{\phi}{\varphi}
 
 \newcommand{\IC}{\Bbb{C}}                     
 \newcommand{\IN}{\Bbb{N}}                     
 \newcommand{\IR}{\Bbb{R}}                     
 \newcommand{\IT}{\Bbb{T}}                     
 \newcommand{\IZ}{\Bbb{Z}}                     

 \theoremstyle{plain} 
 \newtheorem{Theorem}{Theorem}[section]
 \newtheorem{Lemma}[Theorem]{Lemma}

 \newtheorem{Corollary}[Theorem]{Corollary}
 
 \newtheorem{Definition}[Theorem]{Definition}
 
 \theoremstyle{definition} 
  \newtheorem{Remark}[Theorem]{Remark}

 \numberwithin{equation}{section}

\begin{document}

 \title{  Group C*-algebras, Metrics  and an Operator Theoretic Inequality.}
 
 \author{Cristina Antonescu and Erik Christensen}
 \address{Institut for Matematiske Fag University of Copenhagen\\
Universitetsparken 5\\
DK 2100 Copenhagen \O \\
Denmark }
 \email{chris@math.ku.dk,  echris@math.ku.dk}
 \date{\today}


 \begin{abstract}
   
   On a discrete group $G$ a length function may  implement a
   spectral triple on the reduced group C*-algebra. Following A.
   Connes, the Dirac operator of the triple then can induce a metric 
   on the state space of  reduced group C*-algebra. Recent
   studies by M.  Rieffel raise several questions with respect to
   such a metric on the state space. Here it is proven that for a free
   non Abelian group, the metric on the state space is 
   bounded. Further we propose a relaxation in the way a length
   function  is used in the construction of a 
   metric, and we show that for
   groups of rapid decay there are many metrics related to a length
   function which all have all the expected  properties.
   The boundedness result for free groups is based on an estimate of
   the completely bounded norm of a certain Schur multiplier and on
   some techniques concerning free groups due to U. Haagerup. At the end
   we have included a noncommutative version of the Arzel\'{a}-Ascoli
   Theorem.

 \end{abstract}

 \maketitle

 \section{Introduction}
 In the article \cite{Co1} Connes demonstrates that the geodesic
 distance on a compact, spin, Riemannian manifold $\mathcal{M}$ can be
 expressed in terms of an unbounded Fredholm module over the
 C*-algebra ${\rm C}(\mathcal{M})$. The distance between points $p,q$ in
 $\mathcal{M}$ is obtained via the Dirac operator $D$ by the formula
 \begin{displaymath} d(p,q)= \sup \{\vert{a(p)-a(q)} \vert \; | \; a \in
   {\rm C}(\mathcal{M}), \Vert{[D,a]}\Vert\leq {1}\}.
\end{displaymath}
\indent Inspired by the compact manifold $\mathbb{T}$, i.e. the unit
circle, and the well known identity ${\rm C}^{*}(\mathbb{Z})={\rm C}
(\mathbb{T})$,
it is natural to consider discrete groups with length functions and
for such a pair $(G, \ell)$ to define a Dirac operator $D$ on
$l^{2}(G)$ by $(D\xi)(g)=\ell(g)\xi(g).$ This set up is discussed in
\cite{Co1} and Connes proves that $(l^{2}(G), D)$ is an unbounded
Fredholm module for ${\rm C}_{r}^{*}(G)$ if $\ell^{-1}([0,a])$ is a finite set
for each $a \in \mathbb{R}_{+}.$ Still following \cite{Co1}, one can
then try to define a metric on the state
space 
$\mathcal{S}({\rm C}_{r}^{*}(G))$  of ${\rm C}_{r}^{*}(G)$ 
by \begin{displaymath}
  d_{\ell}(\varphi,\psi)=\sup \{\vert{\varphi(a)-\psi(a)} \vert \; |
  \; a \in
  {\rm C}_{r}^{*}(G), \Vert{[D,a]}\Vert\leq {1}\}.
\end{displaymath}
\noindent In the first place it is not clear that
$d_{\ell}(\varphi,\psi)<\infty$ for all pairs, but if that is the case then
one could ask if $\mathcal{S}({\rm C}_{r}^{*}(G))$ is a bounded set 
with respect to this metric, and if so does $d_{\ell}$ generate the
w*-topology on $\mathcal{S}({\rm C}_{r}^{*}(G)).$ 

Especially  Rieffel has studied these questions and some more
general questions concerning the set of all metrics on
$\mathcal{S}({\rm C}_{r}^{*}(G))$, where $G$ is a group with reasonable
properties. We will concentrate on the three questions
mentioned above. We will prove that for a free (non Abelian) group
$\mathbb{F}_{n}$, on finitely many generators $ n$, and $\ell$ the usual length
function on $\mathbb{F}_{n}$, $d_{\ell}$ is bounded. We
think that $d_{\ell}$ induces the w*-topology too, but the
combinatorics involved in proving this seems to be rather
difficult. On the other hand we have realized that if
one allows a variation in the definition of $d_{\ell}$ then one can quite
easily get all the desired properties of $d_{\ell}$. This
change works for the groups of rapid decay \cite{Jo}, a fundamental
concept which will be defined properly in the text. For
a discrete group which is of rapid decay with respect to some
length function $\ell$ there exists a $k\in \mathbb{N}$ such that the
metric $d_{k,\ell}$ on the state space $\mathcal{S}$ of
${\rm C}_{r}^{*}(G)$ defined by
\begin{displaymath} 
d_{k,\ell}(\varphi,\psi)=\sup\{\vert{\varphi(a)-\psi(a)} \vert \; | \;
a \in {\rm C}_{r}^{*}(G),
\Vert{\underbrace{[D,[D,...,[D,a]\dots]}_{k}}\Vert\leq {1}\}.
\end{displaymath}
has all the wanted properties.  The proof of this result is on the
other hand very easy, and this  may be seen as  an indication 
in favor of this
alternative construction of a metric from a length function.

The boundedness
result for $d_{\ell}$ in the case of free groups is based on an
estimate of the completely bounded norm of a certain Schur multiplier.
This is a purely operator theoretic result which is not far from the
known results which are collected in the book by Pisier \cite{Pi}. 
This variant is apparently unknown and a proof of this
result fills section 3. 

At the end of the paper we have
included some thoughts on what a noncommutative version of the
Arzel\'{a}-Ascoli Theorem could be. Using some classical functional
analysis techniques it turns out that at least one way to look at the
Arzel\'{a}-Ascoli Theorem is that if you  know a norm compact convex
balanced subset, say $\mathcal{K}$, of a unital C*-algebra, say
$\mathcal{A}$, such that $\mathcal{K}$ separates the states of
$\mathcal{A}$, then you can tell whether any other subset, say
$\mathcal{H}$ of $\mathcal{A}$ is norm precompact. In order to explain
the result we will use the terminology that
$\mathcal{A}_{\varepsilon}$ denotes the closed ball in $\mathcal{A}$
of radius $\varepsilon$.  Having this, a  subset $\mathcal{H}$
of $\mathcal{A}$ is norm precompact if and only if
\begin{displaymath}
\mathcal{H}\quad  \text{ is bounded and} \quad   \forall \varepsilon >0 \; \exists N\in \mathbb{R}_{+} \quad  \mathcal{H}\subseteq \mathcal{A}_{\varepsilon} 
+N\mathcal{K} + \mathbb{C}I. 
\end{displaymath}
\noindent The result is not deep at all, but it emphasizes that the
concept of a metric on the state space $\mathcal{S}(\mathcal{A})$
generating the w*-topology is closely related to compact convex
subsets of $\mathcal{A}$ and once you know one such set then the  other
ones {\it are not too far away }. This close connection
between metrics on $\mathcal{S}(\mathcal{A})$ and precompact subsets
of $\mathcal{A}$ is studied in details in the papers \cite{Ri2},
\cite{Pav}; from where we have learned it.

 \section{Notation and preliminaries}
 \label{intro}

\medskip

Most of our investigations deal with properties of the C*-algebra
generated by the left regular representation $\lambda$ of a discrete group $G$
on $l^{2}(G)$. We refer to chapter 6 of \cite{KR} for the basic properties of this
C*-algebra, but we will use a slightly different notation which is inspired
by Connes' presentations in \cite{Co1} and \cite{Co2}. 
This means that  for an $x$ in the group algebra $\mathbb{C}G$ 
we will write $\lambda(x)=\sum_{g} x(g)\lambda_{g}$ for the convolution operator
on $l^{2}(G)$, and for $g\in G$, $\delta_{g}$ denotes the natural basis element in $l^{2}(G)$.
The C*-algebra generated by $\lambda(\IC G)$ in $B(l^{2}(G))$ is called
the reduced group C*-algebra and denoted ${\rm C}_{r}^{*}(G)$.
Any element $x$ in this algebra has a unique representation in 
$l^{2}(G)$ by $x\to x\delta_{e}$ so in a natural way we have
\begin{displaymath}
l^{1}(G)\subseteq {\rm C}_{r}^{*}(G)\subseteq l^{2}(G)
\end{displaymath}
and for $x\in l^{1}(G)$
\begin{displaymath}
\lVert{x}\rVert_{2}\le \lVert{\lambda(x)}\rVert \le\lVert{x}\rVert_{1}. 
\end{displaymath}
\indent For a discrete group of rapid decay (RD) one has a kind of
inverse to the first inequality and such a type of inequality is very
powerful as 
we shall see. In order to explain the concept of rapid decay
we remind the reader that a length function $\ell$ on a group $G$ is a
mapping $\ell:G\to\mathbb{R}_{+}\cup \{0\}$ such that

\begin{itemize}
\item[{\it(i)}] $\ell(e)=0$
\item[{\it(ii)}] $\forall g,h \in G: \ell(gh)\le \ell(g)+\ell(h)$
\item[{\it(iii)}] $\forall g \in G: \ell(g)=\ell(g^{-1})$
\end{itemize} 

\begin {Definition}A discrete group G is said to be of rapid decay
  (RD) if there exist
a length function $\ell$ on G and positive reals C,k such that
\begin{displaymath}
\forall x\in \mathbb{C}G: \Vert\lambda(x)\Vert\leq C\Big(\sum_{g\in
  G}(1+\ell(g))^{2k}\vert x_{g}\vert^{2}\Big)^{\frac{1}{2}}.
\end{displaymath}
\end{Definition}
In the very innovative paper \cite{Haa}, Haagerup proved that the free
non Abelian groups 
$\mathbb{F}_{n}$ are of rapid decay with C=2 and k=2.
Also the free Abelian groups $\mathbb{Z}^{k}$ are of rapid decay and
in fact one can get 
an estimate dominating the norm $\Vert{x}\Vert_{1}.$ For k=1 this is
obvious  since
for $x\in \mathbb{C}\mathbb{Z}$
\begin{displaymath}
\sum_{n\in\mathbb{Z}}\vert x_{n}\vert\leq\Big(\sum_{n\in\mathbb{Z}}
(1+\vert n \vert)^{2}\vert x_{n}\vert^{2}\Big)^\frac{1}{2}
\Big(\sum_{n\in\mathbb{Z}}(1+\vert n \vert)^{-2}\Big)^\frac{1}{2}.
\end{displaymath}
The article  \cite{Jo} by  Jolissaint contains a lot of results on
groups of rapid decay and among his results he 
proves that a discrete group is of rapid decay, 
if it  is  of polynomial
growth with respect to some set of generators and the
corresponding  length function. In Connes' book 
\cite{Co2} he presents in 
Theorem 5, p. 241 a proof of the fact that 
 the word hyperbolic groups of Gromov  all are of
rapid decay. 

As mentioned before Connes defines in \cite{Co1} a metric on a
non-commutative C*-algebra  
via an unbounded Fredholm module. For a discrete group $G$ with a
length function $\ell$  he also obtains 
an unbounded Fredholm module if $\ell$ has the sort of finiteness 
property
defined in the next definition. 
\begin{Definition} 
Let $G$  be a discrete group and $\ell:G\to\mathbb{R}_{+}\cup\{0\}$
a length function. If, for each 
$c\in\mathbb{R}_{+}$, $\ell^{-1}([0,c])$ is a finite set then we say
that $\ell$ is proper.  
\end{Definition}
The Fredholm module for ${\rm C}^{*}_{r}(G)$ 
is the Hilbert space $l^{2}(G)$ and 
the Dirac operator $D$ on $l^{2}(G)$ is the selfadjoint unbounded
multiplication  operator which multiplies $\xi\in l^{2}(G)$ by $\ell$ 
 pointwise.
\begin{Definition} \label{dl}
Let G be a discrete group with a length function $\ell$ and let
$\mathcal{S}$ denote the state space  of 
${\rm C}^{*}_{r}(G)$ then $d_{\ell}:\mathcal{S}\times\mathcal{S}\to
[0,\infty]$ is defined by  
\begin{displaymath} 
d_{\ell}(\varphi,\psi)= \sup\{\vert{\varphi(a)-\psi(a)} \vert \; | \;a \in
{\rm C}_{r}^{*}(G), 
\Vert{[D,a]}\Vert\leq {1}\}.
\end{displaymath} 
\end{Definition}
A computation involving the properties of a length function shows that
for any pair  
$g\in G$ and $\xi\in l^{2}(G)$ we have
\begin{displaymath}
\Big([D,\lambda_{g}]\xi\Big)(k)=\Big(\ell(g^{-1}k)-\ell(k)\Big)\xi(g^{-1}k).
\end{displaymath}
So $\Vert[D,\lambda_{g}]\Vert\leq \ell(g)$ and we see that $d_{\ell}$ must
separate the points in $\mathcal{S}.$  It is not clear if
$d_{\ell}(\varphi,\psi)<\infty$ always, but except for that, $d_{\ell}$
behaves exactly as 
a metric on $\mathcal{S}$ so we will call $d_{\ell}$ a possibly infinite
metric on $\mathcal{S}$.

 \section{The completely bounded norm of certain Schur multipliers on $B(H)$}
 \label{cbnorm}
 We will let $M_n(\IC)$ denote the $n \times n$ complex matrices.  A
 matrix $A \in M_n(\IC)$ induces a mapping $A_s : M_n(\IC) \to
 M_n(\IC)$ by the Schur multiplication which is given by the
 expression, $M_n(\IC) \ni X = (x_{ij}) \to  (a_{ij}x_{ij}) = A_s(X) 
 \in M_n(\IC)$. Here we will consider a generalization of this
 mapping to   operators on an abstract Hilbert
 space $H$, which for some discrete group  $G$ 
 is decomposed into a sum of orthogonal subspaces
 $H_g, g \in G, H = \oplus_{g} H_g$. It should be remarked that
 for any 
 $g \in G$ we do permit that $H_g = {0}$ and we want to
 emphasize that the results in this section are designed to
 work for the  well known Abelian group $\IZ$, whereas the work in the
 in rest of the paper aims at general discrete 
groups of rapid decay. Despite the fact that this section really is
devoted to a result on the group $\IZ$  we present the first result in
terms of a general discrete group $G$, a Hilbert space $H$ and a
decomposition of $H$ indexed by $G$. 
In this setting where we have  a decomposition of
 $H = \oplus_{g} H_g$ we get a matrix decomposition of the operators in $B(H)$ such
 that for any $x \; \text{in} \; B(H)$ we can write $x = (x_{s,t}), s,
 t \in G$ where each $ x_{s,t} \; \text{is in} \; B(H_t, H_s) $. It
 is not possible to generalize the Schur multiplication directly to this
 setting since $B(H_t, H_s )$ is not an algebra unless $ s = t $. On
 the other hand it is possible to multiply any operator $x_{s,t} \in B(H_t, H_s)$
 with  a complex scalar, so for an infinite scalar matrix $ \Lambda =
 (\lambda_{s,t}) $ it is possible to perform a formal Schur
 multiplication $ B(\oplus H_g) \ni x = (x_{s,t}) \to 
 (\lambda_{s,t}x_{s,t}) = \Lambda_s(x)$. The latter matrix may not correspond to a
 bounded operator, but  the product is  well defined as an infinite
 matrix and we will call it a formal Schur product. The theorem just
 below, provides a criterion on $\Lambda $ for the boundedness of
 the Schur product $\Lambda_s(x)$ for any $x \; \text{in} \; B(H)$. We
 do a little more since we do compute the so called completely bounded
 norm of $\Lambda_s : B(H) \to B(H)$. The theory connected to
 completely bounded operators is described in Paulsen's book
 \cite{Pau}.
 In order to explain this concept
 shortly  we consider for $n \in \IN$ the mapping $(\Lambda_s)_n : M_n(\IC)
 \otimes B(H) \to M_n(\IC) \otimes B(H)$ which is given as
 $(\Lambda_s)_n = \text{id}_{M_n(\IC)} \otimes \Lambda_s$. If the
 sequence of norms defined by $\| \Lambda_s \|_{n} := \|(\Lambda_s)_n\|$ is
 bounded,  $\Lambda_s$ is
 said to be completely bounded and the completely bounded norm,
 $\|\Lambda_s\|_{cb}$ is given as the $ \sup_n \|\Lambda_s\|_{n}$. Our
 criterion for $\|\Lambda_s\|_{cb}$ to be finite is based
 upon a  generalization of a theorem by  M. Bo\.{z}ejko and
 G. Fendler \cite{BF}. In Pisier's book   \cite{Pi}, he presents this
 result in Theorem 6.4. This theorem 
 deals with the situation where each
 of the summands $H_g$ above is one-dimensional, i. e. $H_g = \IC\;$.

\begin{Theorem} \label{Schur} 
Let $G$ be a discrete group, $H$ be a Hilbert space which is decomposed into a sum of orthogonal
subspaces $H_g, g \in G$ and let $\varphi$ be a complex
function on $G$. If  the 
linear operator $T_{\varphi} : \lambda(\IC G) \to\lambda( \IC G)$
 which is defined by
$ T_{\varphi}(\lambda_g) = \varphi(g)\lambda_g $ extends to a
completely bounded operator  on C$^{\ast}_r(G)$ then for the matrix
$\Lambda$ given by $\Lambda = (\lambda_{s,t})\;  s, t
\in G, \lambda_{s,t} = \varphi(st^{-1}) $ the mapping  $\Lambda_s$ is
completely bounded, $\|\Lambda_s\|_{cb} \leq \|T_{\varphi}\|_{cb}$ and
$\Lambda_s$ is an ultraweakly continuous, or normal, operator on $B(H)$. 

\end{Theorem}

\begin{proof}

Suppose that $T_{\varphi} $ is completely bounded, then for the case 
where  $H_g = \IC$ for every  $g \in G$   the result follows from
\cite{Pi}. The  proof of in that book of the scalar case, as presented in the
proof of the implication [ (i) $ \Rightarrow $ (iii) ] of Theorem 6.4
of \cite{Pi}, shows that there exist a Hilbert 
space K  and two functions say $\xi \; \text{and } \; \eta $ on $G$
with values in $K$ such that 

\begin{displaymath}
  \forall s, t \; \in G \quad  \varphi(st^{-1}) = (\xi(t), \eta(s)) \quad
  \text{and} \quad \forall g \in G \; \|\xi(g)\| \leq \sqrt{\|T_{\varphi}\|_{cb}}\;, \quad \|\eta(g)\| \leq   \sqrt{\|T_{\varphi} \|_{cb}} \;.
\end{displaymath}
We will now turn to the operator $\Lambda_s$ and show that it is
completely bounded. The proof of this can be obtained as a
modification of a part of the proof of Theorem 5.1 of \cite{Pi}. In
fact the representation of $\varphi $, we have obtained above,  makes it
possible to construct operators say $x$ and $y$ in $B(H, H\otimes K)$
such that $\Lambda_s $ can be expressed as a completely bounded
operator in terms of these operators. We recall that $H = \oplus H_g $
and define the operators $x$ and $y$ on a vector $\alpha =
(\alpha_g)_{g \in G}$ by $x\alpha = ( \alpha_g \otimes \xi(g) )_{ g
  \in G}$ and similarly $y\alpha = ( \alpha_g \otimes \eta(g))_{g \in
  G}$ and we find that both operators are of norm at most
$\sqrt{\|T_{\varphi} \|_{cb}}\;.$ Let $\pi$ denote the representation
of $B(H)$ on $H \otimes K$ which is simply the amplification $ a \to a
\otimes I_K$, then an easy computation shows that for any pair of
vectors $ \alpha = ( \alpha_g), \; \beta = ( \beta_g) $ from $H =
\oplus H_g$ and any $a$ in $B(H)$ we have

\begin{align}
  (\pi(a)x\alpha, y\beta) \;&= \; \sum_{s,t \in G}(a \alpha_t, \beta_s)(
  \xi(t), \eta(s))\\
  &= \; \sum_{s,t \in G}(\varphi(st^{-1})a_{s,t}\alpha_t, \beta_s)\\
&= (\Lambda_s(a)\alpha, \beta).
\end{align}
Consequently $\Lambda_s(a) = y^{\ast}\pi(a)x$ and the cb-norm of
$\Lambda_s$ is at most $\| T_{\varphi} \|_{cb}\;.$
The concrete description $\Lambda_s(.) = y^{\ast}\pi(.)x$ where $\pi$
is just an amplification shows that $\Lambda_s$ is ultraeweakly continuous.
\end{proof}

\begin{Corollary}
If $G$ is an Abelian discrete group and $T_{\varphi}$ extends to a
bounded operator on C$^{\ast}_r(G)$ then  the mapping  $\Lambda_s$ is
completely bounded, $\|\Lambda_s\|_{cb} \leq \|T_{\varphi}\|$ and
$\Lambda_s $ is an ultraweakly continuous, or normal, operator on $B(H)$ 
\end{Corollary}

\begin{proof}
We have to prove that $T_{\varphi}$ extends to a completely bounded
mapping if it extends to a bounded mapping and that the two norms on
$T_{\varphi}$ agree. In order to do so we remark that for the compact
Abelian dual
group $\widehat{G}$ and any natural number $k$ we have
C$^{\ast}_r(G)\otimes M_k(\IC) \; =
\;\text{C}(\widehat{G}, M_k(\IC))$, the continuous  $M_k(\IC)$ valued
functions  on $\widehat{G}$.  Then  for any finite sum $ x =
\sum  \lambda_g \otimes m_g $ in $\lambda(\IC G) \otimes M_k(\IC)$ we have
\begin{displaymath}
\|x\| = \max\{ \|\sum_g \chi(g) m_g \|_{M_k(\IC)} \; | \; \chi \in
\widehat{G} \}.
\end{displaymath}
The norm in $M_k(\IC)$ is determined by the
functionals of norm one on this  algebra, so let
$M_k(\IC)^{\ast}_1$ denote this unit ball and we get

\begin{displaymath}
\| x \| = \max\{ |\sum_g \chi( g )\psi( m_g )|  \; | \; \chi \in
\widehat{G} \; \text{and} \; \psi \in M_k(\IC)^{\ast}_1\}.  
\end{displaymath}
Let us now suppose that $T_{\varphi}$ extends to a bounded operator on
the group algebra of norm at most 1. For $x$ as above of norm at most
1 in C$^{\ast}_r(G)\otimes M_k(\IC)$ we get for any pair $\chi, \;
\psi$ as above that $ |\sum_g \chi( g )\psi( m_g )| \leq 1 $, hence
for this fixed $\psi$ we get in C$^{\ast}_r(G)$ the estimate $\|\sum_g
\psi(m_g)\lambda_g\| \leq 1.$ Since $\|T_{\varphi}\|\leq 1 $ we also
get $\|\sum_g \varphi(g) \psi(m_g)\lambda_g\| \leq 1,$ but this holds
for any $\psi$ so we can go back and note that $T_{\varphi}$ is
completely bounded of norm at most 1.

\end{proof}

The following corollary shows that a square summable function
$\varphi$  on a commutative discrete group $G$ induces a completely
bounded operator $\Lambda_s$. Besides these functions and the positive
definite functions  on $G$  we do not know of any other general results which
can guarantee the complete boundedness of $\Lambda_s$.

\begin{Corollary} \label{convolution}

Let G be an Abelian discrete group, $\varphi \in \ell^2(G)$ and let
$\Lambda :G \times G \to \IC $ be given  by $\Lambda(s,t) =
\varphi(st^{-1})$. Then $\Lambda_s$ is completely bounded and
$\|\Lambda_s\|_{cb} \leq \|\varphi\|_2$.

\end{Corollary}
\begin{proof}
The operator $T_{\varphi} $
on C$^{\ast}_r(G)$ can, when we look at the latter algebra as
C$(\widehat{G})$, be expressed as the convolution operator  implemented
 by the Fourier transform,
$\hat{\varphi} \in \text{L}^2(\widehat{G}).$ Here we have chosen the
probability Haar measure on the compact group $\widehat{G}$ such that the
Fourier transform is an isometric  operator between the 2 Hilbert
spaces. It follows directly from The Cauchy-Schwarz inequality that the
norm of 
 the convolution operator is dominated by the norm
 $\|\hat{\varphi}\|_2$, and the corollary follows. 
           
\end{proof}

The purpose of the previous corollary is actually to compute the norm
of the partial inverse of certain derivations on $B(H)$. Let $D$ be a
possibly unbounded self adjoint operator on $B(H)$ with spectrum
contained in the set of integers $\IZ$, then the Hilbert space $H$
decomposes as a direct sum of the eigenspaces say $H_m$ of $D$. Many
of these spaces may vanish, but anyway we can write
$H = \oplus_{ m \in \IZ}H_m$ and we will be able to use the results
from above concerning the norm of certain Schur multipliers. The
question we are going to deal with, is to give a description of  a
bounded operator $a \in B(H)$ which has the property that the
commutator $[ D, a ]$ is bounded and of norm at most 1. Clearly all
bounded operators which commute with $D$ must play a special role in
this set up. This set is a von Neumann algebra and consists of the 
operators in the main diagonal of $B(H)$, when the latter algebra is
viewed as infinite matrices with respect to the decomposition $H =
\oplus H_m $. Consequently we will let $\frak{D_0}$ denote the
commutant of $D$ and for $ k \in \IZ$ we will define the $k$'th
diagonal of $B(H)$ by 

\begin{displaymath}
\frak{D}_k = \{(x_{ij}) \in B(H) \, | \, i - j \neq k \Rightarrow
x_{ij} = 0 \} .
\end{displaymath}
For $k \in \IZ$ there is a natural projection of $B(H)$ onto
$\frak{D}_k$, say $\rho_k$ given by the expressions  

\begin{displaymath}
\rho_k((x_{i,j}))_{m,n} = 
\begin{cases}x_{m,n} & \; \text{if} \; m-n = k \; ,\\
0 & \; \text{if} \; m-n \neq  k \;.\\
\end{cases}
\end{displaymath}
If one computes $\rho_k(x)^{\ast} \rho_k(x)$ it is easy to realize that
$\rho_k$ is a projection onto the $k$'th  diagonal and  of norm at most
one. The problem we are facing is analogous to well known problems
concerning convergence of Fourier series; it is not easy to give norm
estimates of the norm of a general finite sum $\sum_{k \in C}\rho_k(x)$.
If we disregard convergence questions for some time, it follows from
elementary algebraic manipulations that for an operator $a = (a_{m,n})
\in B(H)$, the commutator $[ D, a ]$ must have the formal infinite
matrix $ c = (c_{m,n}) $ given by $c_{m,n} = (m - n )a_{m,n}$. So
at least 
formally we can write 
\begin{displaymath}
[D,a] \; = \; \sum_{k \in \IZ}k\rho_k(a),
\end{displaymath}
and we see that this operator on $B(H)$  in fact is
 an unbounded Schur multiplier. Further
it follows - so far formally - that  
 \begin{displaymath}
a - \rho_0(a)  \; = \; \sum_{k \in \IZ \; \text{and} \;k \neq 0}
\frac{1}{k} \rho_k([ D, a ]),
\end{displaymath}
Hence the partial inverse to the derivation $B(H) \ni a \to [D ,a ]$ is
a Schur multiplier which according to the results above will turn out
to be completely bounded.
With this notation in mind we can offer  norm estimates for such
sums in the next theorem. The theorem is a generalization of the well
known fact that the Fourier series for a differentiable $2\pi$
periodic function on $\IR$ is uniformly convergent. 

\begin{Theorem} \label{commutnorm}
  Let $D$ be a self adjoint operator on a Hilbert space $H$ such that
  the spectrum of $H$ is contained in $\IZ$ and let $K = \{a \in B(H)
\; |\; \| [D,a] \| \leq 1 \; \text{and} \; \rho_0(a) = 0 \}$, then every element
in $K$ is of norm at most $\frac{\pi}{\sqrt{3}}$.  For  $a \in B(H)$ such that  $\| [D,a] \| \leq 1 $ the sum $\sum_{m \in \IZ}\rho_m(a)$ is norm convergent
and $\forall k \in \IN_ :\;\; \|\sum_{ |m| > k}\rho_k(a)\| \leq
\sqrt{\frac{2}{k}}\;.$
 \end{Theorem}

\begin{proof}
  Suppose $a \in B(H)$ satisfies $\| [ D, a] \| \leq 1 \; \text{and}
  \; \rho_0(a) = 0\; $, then the first statement in the theorem is, as
  we shall see, just a special case of the second corresponding to $k
  = 0$, although the estimates are slightly different. Let then $k \in
  \IN_0$ be given and consider the set of Hilbert spaces $H_m, m \in
  \IZ$ where the space $H_m$ is defined as above,
  i. e. the eigenspace for $D$ corresponding to the eigenvalue $m$.
  We can then apply Corollary \ref{convolution} for the group $\IZ$
  and the function $\varphi_k : \IZ \to \IR$ given by

 \begin{displaymath}
\varphi_k(m) = 
\begin{cases}m^{-1}, & \; \text{if} \; |m| > k \; ,\\
0 & \; \text{if} \; |m| \leq  k \; .\\
\end{cases}
\end{displaymath}
Hence we see that 

\begin{displaymath}
\|a - \sum_{i = -k}^k \rho_i(a)\| \leq \| [ D, a ] \| \|\varphi_k\|_2
\leq   (2 \sum_{j > k} j^{-2} )^{\frac{1}{2}}.
\end{displaymath}
After recalling the well known  sum $ \sum_{ j \in \IN} j^{-2} =
\frac{\pi^2}{6}$ for $k=0$ and the integral estimate $ \sum_{ j > k}j^{-2} <
\frac{1}{k}$ for $k>0$, the theorem follows.

\end{proof}
 
\begin{Remark}
The proof of Theorem \ref{commutnorm} depends on the fact that the spectrum of
$D$ is contained in the  Abelian discrete group  $\IZ$. This is not
likely to be a relevant condition for a result of this type and we
would think that this theorem must have a more general version which
is valid for an unbounded  self adjoint operator whose spectrum
consists of points $(s_k)_{k \in \IZ}$ such that $ |s_k| \to \infty  \; \text{
  for }\; |k| \to \infty $ and $\inf\{|s_m - s_n| \; | \; m, \; n \in
\IZ \} > 0.$  We are aware of, and have already
mentioned the fact that an estimate
similar to the one above does
 exist for ordinary differentiation on C$(\IT)$, but
we have not found a general operator theoretic treatment of this
problem. 
\end{Remark} 

\section{On unital C*-algebras as noncommutative compact metric
  spaces.} \label{metric}
This section is mainly devoted to the study of some metrics on the state
space of a 
C*-algebra which is  generated by the left regular representation of a
discrete group of rapid decay.  We will start by recalling Connes'  
construction  \cite{Co1} which defines a metric on the state space of a 
discrete group C*-algebras in
terms of an  unbounded Fredholm module.
Let then  $G$ be a discrete group  with a
 length function $\ell : G \to \IN_0$ such that $\ell$ is
 proper and $G$ is of rapid decay with constants $C, k$ with
 respect to $\ell$. 
 As described in \cite{Co2} p. 241 such a
 length function $\ell$ on a discrete group $G$ with values in $\IN_0
 $ induces a decomposition of $\ell^2(G)$ into an orthogonal sum of
 subspaces $H_m $, each one being the closed linear span of the basis
 vectors $\delta_g$ for which  $\ell(g) = m$. The Dirac operator on
 $\ell^2(G)$ is the self adjoint unbounded operator $D$ which is the
 closure of the operator $D_0$ defined on the linear span of the basis
 vectors $ \delta_g, \; g \in G $ and acts by $D_0\delta_g \; = \;
 \ell(g)\delta_g$.  The operator $D$ clearly has it's spectrum
 contained in $\IZ$ and the eigenspaces  all vanish for $m < 0$ and
 equals $H_m$ for $m \geq 0$. The Theorem \ref{commutnorm} is
 designed to deal with this situation, but it does not work  as well
 as expected. The theorem provides a norm
 estimate of $\|a 
 - \rho_0(a)\|$ in terms of $\|[D,a]\|$ for an $a$ in C$^{\ast}_r(G)$,
 and we believed that from this it would be easy to get an estimate of
 $\|\rho_0(a)\|$, because a group is such a rigid
 object. Unfortunately this is not so and we can only get an estimate
 which also takes care of the main  diagonal part $\rho_0$ if 
the group is a one of the  free non Abelian groups $\mathbb{F}_n$.

\begin{Theorem} \label{dbounded}
  Let $G$ be a free non Abelian 
group on finitely many generators and $\ell$ the
  natural length function on $G$, then the metric $d_{\ell}$ on the
  state space  $\mathcal{S}(${\rm C}$^{\ast}_r(G))$ is bounded, 
and the diameter of the
  state space is at most 5. 
\end{Theorem}

\begin{proof}
We will prove the boundedness of the metric by studying the convex 
and balanced subset
$K$ of   C$^{\ast}_r(G)$ given by  $ K = \{a \in \text{C}^{\ast}_r(G)
\;| \; \|[ D, a ] \| \leq 1 \} $. By Theorem \ref{commutnorm} we know
that for an $a \in K $ we will have $\|a -\rho_0(a)\| \leq
\frac{\pi}{\sqrt{3}} $, so we only have to get an estimate of the norm
$\rho_0(a)$. In order to control
$\|\rho_0(a)\|$ we do first restrict to  the case where $a$ is of finite
support in $G$ and secondly we   make the extra assumption that $(a\delta_e,
\delta_e ) = 0.$  If the latter  is not the case we simply subtract the
corresponding multiple of the unit from $a$. This operation has of
course no effect on the commutator $[ D, a ]$. Since we know by
assumption that this
commutator is of norm at most 1, we get the first estimate 

\begin{displaymath}
 1 \geq  \|[ D, a ]\delta_e \|^2 = \sum_{g \in G}\ell(g)^2|a(g)|^2.
\end{displaymath}
Let us now pick a unit vector $\xi \in H_m$ and let $p_m$ denote the orthogonal
projection from $H$ onto $H_m$. We can now try to estimate $\|\rho_0(a)\|$
by estimating $p_ma\xi$. So let $g  \in G$ be of length $m$ then $g$
can be  expressed uniquely in terms of generators as  $
g = g_1g_2, \dots  ,g_m $, and when we have to compute the value of the
convolution $a*\xi(g)$  we have to remember that $\xi$ is supported on
words of length $m$ so the sum will be an expression of the type 

\begin{displaymath}
a*\xi(g) = \sum_{k=1}^m \sum_{\{s_1, \dots , s_k\;|  s_k \neq g_k \;
  \text{and} \; s_k \neq g_{k+1}^{-1} \} } a(g_1, \dots , g_k s_{k}^{-1}, \dots s_1^{-1})\xi(s_1, \dots , s_kg_{k+1}, \dots g_m),
\end{displaymath}
We can now imitate a trick from the proof of \cite{Haa} Lemma 1.3. In
order to do so we define for a fixed $ k, 1 \leq  k \leq m$ ( remember
$a(e) = 0$ so $ k > 0$ ) 
 a function $b_k$ supported on words of length $ k$ and a vector $\eta_{m-k}$
supported on words of length $m-k$ by 

\begin{align}
b_k(g_1, \dots , g_k) &= (\sum_{\{s_1, \dots , s_k\;|  s_k \neq g_k
  \}}|a(g_1, \dots , g_k s_{k}^{-1}, \dots
s_1^{-1})|^2\;)^{\frac{1}{2}}, \\
b_m(g_1, \dots , g_m) &= |a(g_1, \dots , g_m)| \\
\eta_{m-k}(g_{k+1}, \dots , g_m) &= (\sum_{\{s_1, \dots , s_k\;|  s_k \neq g_{k+1}^{-1}
  \}}|\xi(s_1, \dots , s_kg_{k+1}, \dots g_m)|^2\;)^{\frac{1}{2}} \\
\eta_0(e) &= \|\xi\|_2 = 1.
\end{align}
Having this we get

\begin{displaymath}
|a*\xi(g)| \leq  \sum_{k=1}^m b_k(g_1, \dots ,g_k)\eta_{m-k}(g_{k+1},
  \dots ,g_m) = (\sum_{k=1}^m b_k*\eta_{m-k})(g).
\end{displaymath}
As in \cite{Haa} we will let $\chi_m$ denote the characteristic function
on the words of length $m$ and we will further use the statement
contained in Lemma 1.3  of \cite{Haa} which says that for functions
like $b_k$ which is supported on words of length k and $\eta_{m-k}$
which is supported on words of length $m-k$ one has
\begin{displaymath}
\|b_k*\eta_{m-k}\chi_m\|_2 \leq  \|b_k\|_2 \|\eta_{m-k}\|_2.
\end{displaymath}
A combination of the inequalities above then yields

\begin{align}
\|p_ma\xi\|_2 &= \|a*\xi \chi_m\|_2 \\
              & \leq \|\sum_{k=1}^mb_k*\eta_{m-k} \chi_m\|_2 \\ 
              & \leq \sum_{k=1}^m \|b_k*\eta_{m-k} \chi_m\|_2 \\ 
              & \leq   \sum_{k=1}^m \|b_k\|_2 \|\eta_{m-k}\|_2 \\   
              & \leq   \sum_{k=1}^m \|b_k\|_2, \quad \text{since} \;
              \|\xi\|_2 = 1  \\
              &= \sum_{k=1}^m k\|b_k\|_2(1/k)\\
              &\leq (\sum_{k=1}^m k^2 \|b_k\|_2^2)^{\frac{1}{2}}(\sum_{k=1}^m k^{-2})^{\frac{1}{2}}\\
              &\leq \frac{\pi}{\sqrt{6}} (\sum_{k=1}^m k^2 \sum_{g,
                \ell(g) = 2k} |a(g)|^2)^{\frac{1}{2}}\\
              & \leq \frac{\pi}{\sqrt{6}} (\frac{1}{4} \sum_{g \in
                G}\ell(g)^2 |a(g)|^2 )^{\frac{1}{2}}\\
              & \leq \frac{\pi}{2\sqrt6}.
\end{align}
The computations just above shows that for an $a \in K$ we have $\|\rho_0(a)
-(a\delta_e,\delta_e)I\| \leq \frac{\pi}{2\sqrt6}$ and from
 Theorem \ref{commutnorm} we know that $\| a - \rho_0(a)\| \leq
 \frac{\pi}{\sqrt3}$, so we obtain
$\|a -(a\delta_e,\delta_e) I \| < 2.5$, and  the diameter of the
state space is at most 5.

\end{proof}

In \cite{Ri2}, Rieffel introduces the concept of a lower
semicontinuous Lipschitz seminorm $L$ on a C*-algebra $\mathcal{A}$. The term
{\em Lipschitz} means that the kernel of the seminorm consists of the scalars
and the terms {\em lower semicontinuous } means that the set $\{ a \in
\mathcal{A}\; | \; L(a) \leq 1 \} $ is norm closed.
In our context the operator $D$ induces several   
Lipschitz seminorms 
whose domains of definition always contain the dense subalgebra
$\lambda(\IC G)$  of C$^{\ast}_r(G)$.
In order to define these seminorms we fix the setting as above. 
Let $G$ be a discrete group with a length function $\ell :G \to
  \IN_0$ such that $\ell^{-1}(0) = \{e\}$. Let $D$ denote the
corresponding Dirac operator and $\delta$ the unbounded derivation on
C$^{\ast}_r(G)$ given by $\delta(a) = \rm{closure}([D, a])$, 
if the  commutator 
$ [D,a] $ is bounded and densely defined.  For any natural number $k$
we define a seminorm $L_D^k$ by 
\begin{displaymath}
\text{domain}(L_D^k) \; = \; \text{domain}(\delta^k) \quad \text{and} \quad
 L_D^k(a) = \|\delta^k(a)\|.
\end{displaymath}
 Having this notation we can state a  theorem of quite general
validity.

\begin{Theorem} \label{lsc}
  Let $G$ be a discrete group with a length function $\ell :G \to
  \IN_0$ such that $\ell^{-1}(0) = \{e\}$.  For any natural number $k$
  the seminorm $L_D^k$ on C$^{\ast}_r(G)$ is a lower
  semicontinuous Lipschitz seminorm.
\end{Theorem}

\begin{proof}
The condition  $\ell^{-1}(0) = \{e\}$ implies that the only operators
in  C$^{\ast}_r(G)$ which commute with $D$ are the multiples of the unit
in  C$^{\ast}_r(G)$. Let us define
$\varphi:  \IZ \to \IR $ by

 \begin{displaymath}
\varphi(m) = 
\begin{cases}m^{-1}, & \; \text{if} \; m \neq 0 \; ,\\
0 & \; \text{if} \; m = 0 \; .\\
\end{cases}
\end{displaymath}
 and let $\Lambda_s$ denote the Schur multiplier on $B(\oplus H_m)$
implemented by the function $\lambda(m, n) = \varphi(m-n)$. Then
by Theorem \ref{commutnorm} we know that $\Lambda_s$ is a completely 
bounded
and ultraweakly continuous operator on $B(H)$.
 Let now  $B(H)_1 $ denote the unit ball in
$B(H)$ then for any $k \in \IN$ we have $\Lambda_s^k(B(H)_1)$ is
ultraweakly compact. We can now control most of the set  
\begin{displaymath}
\{ a \in C^{\ast}_{r}(G) \; | \; L_D^k(a) \leq 1\}
\end{displaymath}
the only part missing is the main diagonal $
\rho_0({\rm C}^{\ast}_{r}(G))$.
For $B(H)$ we have $\rho_0(B(H)) = \mathcal{D}_0 \; $ i. e. the main
diagonal which is clearly ultraweakly closed. The sum
$\Lambda_s^k(B(H)_1) + \rho_0(B(H))$ is then ultraweakly closed and
consequently also norm closed and the intersection below is  norm
closed too.

\begin{displaymath} 
\{ a \in \; \text{C}^{\ast}_r(G) \; | \;L^k_D(a)  \leq 1 \} =
[\Lambda_s^k(B(H)_1) + \rho_{0}(B(H))]\cap\text{C}^{\ast}_r(G). 
\end{displaymath}
\end{proof}

It is rather easy to check that the metric $d_{\ell}$ introduced in
Definition \ref{dl} induces a topology which is finer than the
w*-topology. In fact the norm dense group algebra $\lambda(\IC G)$ is
obviously contained in the domain of definition for the derivation
$\delta$ and the metric clearly induces a topology on the state space
which is finer than pointwise convergence on the operators
$\lambda_g$.  The question is whether the 2 topologies agree. In the
first place we thought that once the boundedness question was settled
this question ought to be easy to settle because Theorem
\ref{commutnorm} controls problems involving norms of the diagonals
with large indices. A closer analysis shows that the problems involved
seems to be much more complex and probably are of a difficult
combinatorial nature. The short formulation of the problem is that for
a group element - say $g$ - such that $\ell(g) $ is ``large'' the
unitary operator $\lambda_g$ may have a lot of non vanishing diagonals
with ``small'' indices. At least in principle this makes it possible
for the algebra $\IC G$ to have the property that for any natural
number $N$, there exists an operator $x = \sum x(g)\lambda_g$ such
that $x(g) \neq 0 \Rightarrow \ell(g) > N$, and $\| \sum_{ |k| \leq N
  } \rho_k(x)\| $ is {\em big } whereas $\| \sum_{ |k| > N }
\rho_k(x)\| $ is {\em small}. Since we have not been able to solve
these problems we have looked for alternative constructions of metrics
which will induce the w*-topology on the state space of
C$^{\ast}_r(G)$. The most obvious thing to do, seemed to be to restrict
the attention to the analysis of discrete groups of rapid decay
\cite{Co2},\cite{Jo}.

Before we state and prove our  result
we want to mention that it follows from Rieffel's works \cite {Ri1}, \cite{Ri2}
and the work of Pavlovi\'{c} \cite{Pav} that a lower semicontinuous 
 Lipschitz seminorm L on a unital C*-algebra 
$\mathcal{A}$ is bounded and induces the w*-topology on the state space of $\mathcal{A}$
if and only if the set
\begin{displaymath}
\{a\in \mathcal{A}: L(a) \leq 1\}
\end{displaymath}
has a compact image in the quotient space $\mathcal{A}/{\mathbb{C}I}$,
equipped with the quotient norm.

\begin{Theorem}
  Let G be a discrete group with a length function $\ell :G\to \IN_0$
  such that $\ell^{-1}(e)=0, \; \ell$ is proper and $G$ is of
  rapid decay with respect to $\ell$. Then there exists a
  $k_{0}\in\mathbb{N}$ such that for all $k \geq k_{0}, \; L^k_D $ is
  lower semicontinuous, the metric generated by $L^{k}_D$ on
  $\mathcal{S}({\rm C}^{*}_{r}(G))$ is bounded and the topology
  generated by the metric equals the w*-topology.
\end{Theorem}
\begin{proof}
  
  The statement on lower semi continuity follows from Theorem \ref{lsc}
  and the assumption of rapid decay implies that there exist two
  positive reals $C, s$ such that

\begin{displaymath}
\forall x\in\mathbb{C}G \quad
\Vert\lambda(x)\Vert\leq C\big(\sum_{g}\big( 1+l(g)\big)^{2s}\vert\lambda(g)\vert^{2}\big)^\frac{1}{2}
\end{displaymath}
The number $k_0$ is then defined by  $k_0 = \lfloor s\rfloor +1$, and
given this we will fix a $k \in \IN$ such that $k \geq k_0$.
According to the statement just in front of this theorem we have to
prove that the set, say $\tilde{\mathcal{K}}_k$ defined by  
\begin{displaymath}
\tilde{\mathcal{K}}_k = \{a\in {\rm C}^{*}_{r}(G)\; | \; L^k_D(a) \leq 1\}
\end{displaymath}
has precompact image in ${\rm C}^{*}_{r}(G)/\IC I$. The way we obtain
this is by choosing the element from each equivalence class in
$\tilde{\mathcal{K}}_k$ which is of trace 0. This set is denoted $\mathcal{K}_k$ and is clearly precompact if and
only if $\tilde{\mathcal{K}}_k$ has precompact image in the quotient space
${\rm C}^{*}_{r}(G)/\IC I$. Consequently  $\mathcal{K}_k $ is given by 
\begin{displaymath}
\mathcal{K}_k = \{a\in {\rm C}^{*}_{r}(G)\; | \; L^k_D(a) \leq 1 \quad
\rm{and} \quad (a\delta_e,\delta_e)= 0\}.
\end{displaymath}
The first observation we need has already been used before, namely
that any element $a \in {\rm C}^{*}_{r}(G)$ can be expressed as an $l^2$
convergent infinite sum $\sum a(g)\delta_g$ and that $\|a\|_2 =
\|a\delta_e\|$. Having this, and the fact that $D\delta_e = 0$
 we get for an $a \in \mathcal{K}_k$ that 
\begin{displaymath}
 1 \geq L^k_D(a) =
 \|\delta^k(a)\| \geq \|\delta^k(a)\delta_e\| = \| \sum
 \ell(g)^ka(g)\delta_g \|.
\end{displaymath}
 In particular we get for an  $a \in
 \mathcal{K}_k$ that 
\begin{displaymath}
\sum \ell(g)^{2k} |a(g)|^{2}\leq 1.
\end{displaymath}
The properness condition on $\ell(g)$ implies that there are
only finitely many group elements of length less than any natural
number $n$. Hence in order to prove that $\mathcal{K}_k$ is precompact
it is sufficient to show that for any positive real $\varepsilon$
there exists a natural number $n$ such that for any $a \in
\mathcal{K}_k$ 
\begin{displaymath}
\|\sum_{\ell(g) \geq n} a(g)\lambda_g\|_{{\rm C}^{*}_{r}(G)} \leq \varepsilon
\end{displaymath}
but this is on the other hand easily obtainable from the inequality at
the top of the proof.  In fact let $n \in \IN$ then for $ g \in G \;
\text{with} \; \ell(g) \geq n \geq 1$ we get  
\begin{displaymath}
(1 + \ell(g))^{2s} \leq
2^{2s}\ell(g)^{2s} \leq 2^{2s}n^{(2s - 2k)}\ell(g)^{2k}.
\end{displaymath}
 Since $2s-2k <
0$ there exists an $n \in \IN$ such that $2^{2s}n^{(2s - 2k)} \leq
\frac{\varepsilon^2}{C^2}$.  For this $n $ we then obtain
\begin{displaymath}
\|\sum_{\ell(g) \geq n} a(g)\lambda_g\|_{{\rm C}^{*}_{r}(G)}^2  \leq
C^2\sum_{\ell(g) \geq n}\big(1+\ell(g)\big)^{2s}|a(g)|^2 \leq
C^2\frac{\varepsilon^2}{C^2}L^k_D(a)^2  \leq \varepsilon^2 
\end{displaymath}
and the theorem follows. 
\end{proof}

Given the common use of language which says:  {\em a noncommutative compact
  topological space is a unital C*-algebra}, it seems  natural to
propose a definition of a metric on this noncommutative 
space in terms of an object
which is related directly to the algebra and not to its state space.  
If there is an obvious smooth structure in
terms of a spectral triple, this object should be preferred
since it contains much more information,
 but  for a general noncommutative, unital
and separable 
C*-algebra $\mathcal{A}$ 
without any particularities it seems that any precompact  balanced and
convex subset of $\mathcal{A}$ 
which separates the states on $\mathcal{A}$ contains all the
information needed. The reason why we propose
to look at such a set as a noncommutative metric 
 is partly due to  our reading of  the works by Rieffel and Pavlovi\'{c},
 but also
 partly due to the result in our next section. 

\begin{Definition}
Let $\mathcal{A}$ be a unital C*-algebra.  A subset
$\mathcal{K}$ of  $\mathcal{A}$ is called a metric set if it is norm
compact, balanced, convex  and separates the states on  $\mathcal{A}$.
\end{Definition}
The term balanced will be defined properly in the next section. 
With this definition at hand  one can easily construct metric sets for
separable unital C*-algebras and for instance for a countable group 
$G = \{g_n \; | \; n \in \IN \}$ a metric set in C$^{\ast}_r(G)$ 
could be given by the following expression where $\overline{{\rm
    conv}}$ means the closed convex hull.

\begin{displaymath}
\mathcal{K} := \overline{{\rm conv}}\big(\cup_{n = 1}^{\infty} \{\alpha
\lambda_{g_n} + \beta \lambda_{g_n}^{\ast}\; |
\; \alpha, \beta \in \IC \quad {\rm and} \quad |\alpha| + |\beta| \leq 1/n \}\big).
\end{displaymath}
The next section will contain proofs which hopefully will justify this
introduction of yet another concept.

\section{A noncommutative Arzel\'{a}-Ascoli Theorem}
\label{sec:ncascoli}

The classical Arzel\'{a}-Ascoli Theorem gives a characterization of
precompact subsets of ${\rm C}(X)$ for a compact topological space $X$.
If $X$ is a equipped with a metric $\rho$ generating the topology on
$X$ one can construct a convex subset $\tilde{\mathcal{K}}$ of  ${\rm C}(X)$ by

\begin{displaymath} 
  \tilde{\mathcal{K}} = \{ f \in {\rm C}(X) \; | \; \forall x,y \; \in
  X \;\;
  |f(x) - f(y)| \leq \rho(x,y) \; \}
\end{displaymath}

This set is unbounded since any constant function belongs to
$\tilde{\mathcal{K}}$.  If one normalizes the set by considering the
subset consisting of those elements which all vanish at a certain
point $x_0$ then the classical Arzel\'{a}-Ascoli Theorem shows that the
set, say  $\mathcal{K}$,  given by 
\begin{displaymath} 
  \mathcal{K} = \{ f \in {\rm C}(X) \; | \; \forall x,y \; \in X \;\;
  |f(x) - f(y)| \leq \rho(x,y) \; {\rm and} \; f(x_0) = 0 \}
\end{displaymath}
will be a compact balanced convex subset of ${\rm C}(X)$ which
separates the points in $X$. The Arzel\'{a}-Ascoli-Theorem measures any
other subset of ${\rm C}(X)$ against this set in order to see whether
this subset is precompact or not. What we do in the following lines is
just to transfer this measuring process to the noncommutative case.
The methods we use are elementary functional analytic duality results.
So we have wondered if this sort of result is valid in a much wider
generality like operator spaces \cite{ER}, \cite{Ke}. 
It seems that the validity of a generalization of Lemma 5.1 to this
new setting is the
crucial thing.   Before we start we want to introduce some more
notation. We will be considering the self adjoint part
 of a unital C*-algebra which is denoted $\mathcal{A}_h$
and we want to think of the elements in $\mathcal{A}$ as  
 affine complex 
w*-continuous functions on the state space $\mathcal{S}$ of
$\mathcal{A}$, so we will let $\Bbb{A}(\mathcal{S})$ denote 
the space of w*-continuous affine complex  functions on $\mathcal{S}$
and for an element $a \in \mathcal{A}, \; \widehat{a}$
 will denote the
corresponding affine function in  $\Bbb{A}(\mathcal{S})$. This
presentation of $\mathcal{A}$ is called Kadison's functional
representation of $\mathcal{A}$. It is well known that the functional
representation is isometric on $\mathcal{A}_h$, but for a  general
element $ a \in \mathcal{A}$ we only have the estimates 

\begin{displaymath}
 \|a\| \geq 
\sup|\widehat{a}(\varphi)| = \|\widehat{a}\| \geq \frac{1}{2}\|a\|.
\end{displaymath}
In particular this shows that a subset $\mathcal{H}$ of $\mathcal{A}$
is bounded if and only if the subset $\widehat{\mathcal{H}}$ of 
$\Bbb{A}(\mathcal{S})$ is bounded.
 
 The term
balanced is used in the sense that a subset $\mathcal{H}$ of
$\mathcal{A}$ is balanced if for any complex number $\mu $ such that 
$| \mu | \leq 1 $ we have $\mu \mathcal{H} \subseteq \mathcal{H}$.
 We
remind the reader that a subscript attached to a Banach space like
$Y_{\mu}$ means that we consider the closed ball of radius $\mu$ in
$Y$ and an asterix  attached to $Y$ like $Y^{\ast}$ means the dual
space. For a pair  of Banach spaces like $\mathcal{A}$ and
$\mathcal{A}^{\ast}$ we will use the duality result known under the
name of the bipolar theorem. Here the polar of a set $\mathcal{H}
\subseteq  \mathcal{A}$ is
denoted  $\mathcal{H}^{\circ}$ and defined by  

\begin{displaymath}
\mathcal{H}^{\circ} = \{ \gamma \in \mathcal{A}^{\ast} \; | \; \forall
h \in \mathcal{H} \; |\gamma(h)| \leq 1 \; \}
\end{displaymath}
The bipolar theorem with respect to this polar,
 then states that  the bipolar $\mathcal{H}^{\circ \circ}$,
which now is a subset of  $\mathcal{A}$, is the smallest balanced,
convex and norm closed set in  $\mathcal{A}$ which contains  $\mathcal{H}$.

We can now start the presentation of the generalization of the
Arzel\'{a}-Ascoli Theorem 
 and our first lemma  is closely connected
to the  very fundamental structure in C*-algebra theory, that a
continuous real linear functional can be decomposed in a unique way
into a difference of two positive functionals, such that a certain norm
identity holds.

 \begin{Lemma}
Let $\mathcal{A}$ be a unital  C*-algebra and $\mathcal{S}$ the state space of $\mathcal{A}$, then
\begin{displaymath}
\mathcal{S}-\mathcal{S}=(\mathcal{A}_{h}^{*})_{2}\cap\{\IC I\}^{\bot}.
\end{displaymath}
\end{Lemma}
\begin{proof}
  The inclusion "$\subseteq$" is obvious. To prove the remaining
  inclusion "$\supseteq$" let us take an arbitrary element $f$ in
  $(\mathcal{A}_{h}^{*})_{2}\cap\{\IC I \}^{\bot}$. It is well
  known that for  $f$  in $(\mathcal{A}_{h}^{*})_{2}$ we can find two
  positive linear functionals $f^{+}$ and $f^{-}$ such that
  $f=f^{+}-f^{-}$ and $\Vert f\Vert=\Vert f^{+}\Vert +\Vert
  f^{-}\Vert$. If $f = 0$ we can write $f$ as a difference $ g- g$
  where g is any state on the unital algebra $\mathcal{A}$. 
If $f \neq 0$ the condition $f(I) = 0$ implies that $0 \neq \|f^{+}\| =
\|f^{-}\| = \frac{1}{2}\|f\| \leq 1$. Based on $f^{+}$ we can then
define a positive functional $g$ of norm $\|g\| = 1 - \|f^{+}\|$ by
\begin{displaymath}
g=\frac{(1-\Vert f^{+}\Vert)}{\Vert f^{+}\Vert}f^{+}.
\end{displaymath} 
By construction it follows that   $f^{+}+g$ and $f^{-}+ g$ are
 both states and  from the equality
\begin{displaymath}
f=(f^{+}+g)-(f^{-}+g)
\end{displaymath}
we can conclude that $f \in \mathcal{S}-\mathcal{S}$, and the 
lemma follows.
\end{proof}

Following the ideas of Rieffel \cite{Ri1} we deduce the following result
\begin{Lemma}
Let $\mathcal{A}$ be a unital C*-algebra, $\mathcal{S}$ the state
space of $\mathcal{A}$ and $\mathcal{K}$
a norm  compact subset of $\mathcal{A}$ which separates
the points in the state space. Then for states $\varphi, \psi $ on
$\mathcal{A}$ the  formula 
\begin{displaymath}
d_{\mathcal{K}}(\varphi,\psi):= \sup_{k\in\mathcal{K}}\vert(\varphi-\psi)(k)\vert
\end{displaymath}
defines a metric on the state space 
$\mathcal{S}$ which generates the w*-topology. 
\end{Lemma}
\begin{proof}
The separation property  and the compactness assumption show that 
 $d_{\mathcal{K}}$ is a bounded metric on
$\mathcal{S}$. The norm compactness of $ \mathcal{K}$ and the
boundedness of $\mathcal{S}$ further implies that the topology 
induced by  
$d_{\mathcal{K}}$ is a Hausdorff topology  weaker
 than the compact w*-topology on
  $\mathcal{S}$. A well known theorem from topology then tells 
that the  two topologies do
  agree, and the lemma follows.  
\end{proof}

We can now state and prove the result of this section. 

\begin{Theorem}
Let $\mathcal{A}$ be a unital C*-algebra and $\mathcal{K}$
a  metric subset of $\mathcal{A}.$ 
For any subset $\mathcal{H}$ of $\mathcal{A}$
 the following conditions are equivalent
\begin{itemize}

\item[{\it(i)}] The set $\mathcal{H}$ is norm precompact.

\item[{\it(ii)}] The set of affine functions 
$\{ \widehat{h} \in \Bbb{A}(\mathcal{S}) \; | \; h \in \mathcal{H} \}$
is bounded and equicontinuous with respect to the w*-topology on $\mathcal{S}$.

\item[{\it(iii)}] The set $\mathcal{H}$ is bounded and for 
every $\varepsilon >0$ there exists a real $N>0$ such that
\begin{displaymath}
\mathcal{H}\subseteq\mathcal{A}_{\varepsilon}+N\mathcal{K}+\IC I.
\end{displaymath}
\end{itemize}
\end{Theorem}
\begin{proof}

The equivalence between $(i)$ and $(ii)$ follows from the classical
Arzel\'{a}-Ascoli Theorem 
and the fact mentioned above that $\mathcal{H}$ is bounded if and only
if $\widehat{\mathcal{H}}$ is bounded. 
 
To prove the equivalence between $(ii)$ and $(iii)$ we start with 
$(iii)\Rightarrow (ii)$. From the boundedness of 
$\mathcal{H}$ it follows
that the set $\{ \widehat{h} \; | \; h \in \mathcal{H}
\}$ is bounded. To prove the equicontinuity of this set  
let us fix an $\varepsilon >0$ and find a positive real  $N$ which
fulfills the condition $(iii)$ with respect to $\frac{\varepsilon}{4}$.
Moreover let  $\varphi$ and $\psi$ be two states such that 
\begin{displaymath}
d_{\mathcal{K}}(\varphi,\psi)  < \frac{\varepsilon}{2N}.
\end{displaymath}
then we will show that for any $h \in \mathcal{H}, \;
|\widehat{h}(\varphi) - \widehat{h}(\psi)| \leq \varepsilon$.
Let now   $h$ be  an arbitrary element in 
$\mathcal{H}$, then by $(iii)$ we can find an element
$a\in\mathcal{A}_{1}$,  an element
$k\in\mathcal{K}$ and a complex number $\mu$ such that
\begin{displaymath}
h=\frac{\varepsilon}{4}a+Nk+\mu 1.
\end{displaymath}
then we  obtain
\begin{align}
\vert \widehat{h}(\varphi)-\widehat{h}(\psi)\vert&=\vert
(\varphi-\psi)(h)\vert \\
&\leq
\vert (\varphi-\psi)(\frac{\varepsilon}{4}a)\vert +\vert
(\varphi-\psi)(Nk)\vert \\
&\leq \frac{\varepsilon}{2} +Nd_{\mathcal{K}}(\varphi,\psi) \\
&<\varepsilon.
\end{align}
and the  equicontinuity of $\widehat{\mathcal{H}}$ has been established.

To prove the last implication  $(ii)\Rightarrow
(iii)$ we again first mention that $\mathcal{H}$ is bounded. 
Let then   $\varepsilon >0$ be given and find, via the equicontinuity
assumption on $\widehat{\mathcal{H}}$, a  $\delta >0$ such that
\begin{displaymath}
\forall h \in \mathcal{H} \; \forall \varphi, \psi \in \mathcal{S}: \quad
d_{\mathcal{K}}(\varphi,\psi)\leq\delta \Rightarrow \vert
(\varphi-\psi)(h)\vert = 
\vert \widehat{h}(\varphi)-\widehat{h}(\psi)\vert 
 \leq \varepsilon.
\end{displaymath}
We will now use the bipolar theorem  and remark that the expression 
$d_{\mathcal{K}}(\varphi,\psi)\leq\delta $ exactly means that $\varphi
- \psi \in \delta(K^{\circ})$ It is clear that $\varphi - \psi \in
\mathcal{S} - \mathcal{S}$ and an application of Lemma 5.1 then shows
that the implication above can just as well be expressed as

\begin{displaymath}
\forall h \in \mathcal{H} \;  \forall \gamma \in (\mathcal{A}^{\ast}_h)_2
\cap\{\IC I\}^{\bot}\cap\delta(\mathcal{K}^{\circ}) : \;\;
|\gamma(h)| \leq \varepsilon.
\end{displaymath}
This statement is not sufficient for our computations because it
involves the space $\mathcal{A}_h^{\ast}$ rather that just  
$\mathcal{A^{\ast}}$. Since a functional on $\mathcal{A}$ vanishes on
the identity $I$ if and only both it's hermitian and it's skew hermitian
part vanish on $I$, we can change from $\mathcal{A}_h^{\ast}$ to
$\mathcal{A}^{\ast}$ at the compensation of a factor of 2, so we have
\begin{displaymath}
\forall h \in \mathcal{H} \;  \forall \gamma \in \mathcal{A}_{2}^{*}
\cap\{\IC I\}^{\bot}\cap\delta(\mathcal{K}^{\circ}) : \;\;
|\gamma(h)| \leq 2\varepsilon.
\end{displaymath}
Since all the sets involved  now are convex and  balanced
the bipolar theorem can be applied very easily. Moreover
$\mathcal{K}$ is norm compact so any set of the form
$\mathcal{A}_{\varepsilon}+\IC I+N\mathcal{K}$ is norm closed balanced
and convex. The relation we just established gives immediately the
first inclusion below and the rest follows by some well known ``polar
techniques'' and an application of the bipolar theorem.

\begin{align}
\mathcal{H}& \subseteq 2\varepsilon\Big(\mathcal{A}_2^{\ast}\cap\{\IC I\}
^{\bot}\cap\delta(\mathcal{K}^{\circ})\Big)
^{\circ} \\  
&= 2\varepsilon\Big(\mathcal{A}_{\frac{1}{2}}\cup\ \IC I
\cup \frac{1}{\delta} \mathcal{K}\Big)^{\circ \circ} \\
&= 2\varepsilon\overline{\rm{conv}}\Big(
\mathcal{A}_{\frac{1}{2}}
\cup \IC I \cup\frac{1}{\delta}\mathcal{K} \Big)\\
&\subseteq 2\varepsilon\Big(\mathcal{A}_{\frac{1}{2}}+\IC I+
\frac{1}{\delta}\mathcal{K}\Big).
\end{align}
In conclusion for the given $\varepsilon$ we found 
the number $N=\frac{2\varepsilon}{\delta}$
such that
\begin{displaymath}
\mathcal{H}\subseteq\mathcal{A}_{\varepsilon}+\IC I+N\mathcal{K}
\end{displaymath}
which proves the desired implication.
\end{proof}

\end{document}